\RequirePackage{lineno}
\documentclass[aps,prb,preprint,groupedaddress,showkeys]{revtex4}
\usepackage[pdftex]{graphicx}
\usepackage{amsmath,amssymb}
\usepackage{booktabs}
\usepackage{url}
\usepackage[colorlinks=true,linkcolor=blue,citecolor=blue,urlcolor=blue]{hyperref}
\graphicspath{{figures/}}
\DeclareGraphicsExtensions{.pdf,.jpeg,.jpg,.png}

\begin{document}

\def\papertitle{Quantum-Inspired Hamiltonian Descent for Radiation Therapy Treatment Planning}
\title{\papertitle}
\author{Mohammadhossein Mohammadisiahroydi}
\email[Email: ]{mhms379@umbc.edu}
\thanks{Author to whom correspondence should be addressed.}
\affiliation{University of Maryland, Baltimore County, MD, USA}
\author{Yuxiang Peng}
\affiliation{Purdue University, IN, USA}
\author{Jiaqi Leng}
\affiliation{Viginia Tech, VA, USA}
\author{Robabeh Rahimi}
\affiliation{Inova Health System, VA, USA}
\author{Arezoo Modiri}
\affiliation{University of Maryland, School of Medicine, MD, USA}
\author{Wei Zou}
\affiliation{University of Pennsylvania, Perelman School of Medicine, PA, USA}
\author{Tam\'as Terlaky}
\affiliation{Lehigh University, PA, USA}
\author{Xiaodi Wu}
\affiliation{University of Maryland, College Park, MD, USA}
\date{\today}

\begin{abstract}
\textbf{Background:} Direct incorporation of dose--volume constraints into inverse radiotherapy planning leads to large optimization formulations that faithfully represent clinical objectives but computationally burden the clinical workflow. In current approaches, finding the optimal solution is  compromised by multiple levels of simplification, relaxation and approximation. 

\textbf{Purpose:} To develop and computationally evaluate a quantum-inspired approach for automated dose--volume-constrained intensity-modulated radiation therapy (IMRT) and intensity-modulated proton therapy (IMPT) planning. 

\textbf{Methods:} Prescription and protocol requirements were formulated as mixed integer optimization problems, with binary variables representing potential violations of dose volume constraints. Patient specific dose influence matrices were generated for photon therapy at the control point level, and for proton therapy at proton spot level. The corresponding  intensities were optimized using Quantum Inspired Hamiltonian Descent (QIHD). The proposed approach was evaluated on eight clinical cases, including four IMRT cases (two prostate, one head and neck, and one liver) and four IMPT cases (two lung, one prostate, and one liver). Computational performance was assessed through problem size, solution time, and dose volume histograms (DVHs). 

\textbf{Results:} The IMRT formulations contained 27,621--112,565 variables and were solved by QIHD in 5--32 min on an NVIDIA RTX 4060 GPU. The IMPT formulations contained 69,861--221,342 variables and were solved in 11--35 min. QIHD generated dose distributions that satisfied the clinical protocols within a shorter processing time compared to the current clinical process and with no manual intervention (customary in current treatment planning). QIHD demonstrated applicability for both photon and proton planning. A CPU implementation of Gurobi, a widely used optimization platform, failed to achieve a solution for the same mixed-integer formulations within several hours for most of our tested patient cases.

\textbf{Conclusions:} GPU-accelerated QIHD provides a scalable computational approach for solving dose--volume-constrained mixed-integer treatment-planning formulations and supports an automated workflow for both IMRT and IMPT. The present results establish computational feasibility on clinically sized mathematical formulations. Prospective clinical validation and execution of decomposed subproblems on quantum hardware remain the subject of future work.
\end{abstract}

\maketitle

\section{Introduction}
Directly encoding dose--volume constraints for radiotherapy planning leads to mixed-integer optimization (MIO) formulations where binary variables indicate whether individual voxels exceed prescribed thresholds. Such formulations represent the clinical requirements explicitly but are computationally challenging at realistic scale. \cite{zinchenko2008controlling,shepard1999optimizing} Modern intensity-modulated techniques present such a high-dimensional complex optimization problem. \cite{mohammadisiahroudi2024optimization, mohammadisiahroudiadaptive} 

Clinical treatment planning therefore commonly combines tractable surrogate objectives with iterative adjustment of priorities and constraints. Automated and multicriteria planning methods have reduced manual effort in several settings, but the relationship between the optimized surrogate objective and the final dose--volume requirements remains an important practical consideration.\cite{craft2012multicriteria,breedveld2012icycle,heijmen2018fully,hussein2018automation,mcintosh2017fully} The computational burden can be even greater for intensity-modulated proton therapy (IMPT), where fine spatial dose deposition and robustness considerations produce large optimization problems.\cite{bortfeld2008robust,unkelbach2018robust}

Quantum computing offers a promising approach to addressing the
combinatorial complexity of radiotherapy treatment planning,
particularly for optimization problems involving discrete decisions
and dose--volume constraints.\cite{rahimi2026quantum}
Recent studies have also explored quantum-inspired optimization
algorithms for radiotherapy planning, demonstrating the potential
of alternative computational approaches to these challenging
problems.\cite{pakela2020quantum}

Among emerging quantum optimization approaches, Quantum Hamiltonian
Descent (QHD) provides a Hamiltonian-dynamics framework for
continuous optimization \cite{leng2023quantum,kushnir2025qhdopt}. Current quantum devices, however, cannot directly represent the tens to hundreds of thousands of variables occurring in clinical-scale dose optimization. We therefore investigate Quantum-Inspired Hamiltonian Descent (QIHD), a classical Hamiltonian-dynamics algorithm motivated by QHD and designed for massively parallel GPU execution.\cite{chaudhary2025quantum} QIHD is not a quantum computation and the experiments reported here do not claim quantum speedup. Instead, it provides a scalable method for solving full-size treatment-planning models while preserving a methodological connection to a future hybrid QHD workflow.

The purpose of this study was to (1) formulate an automated dose--volume-constrained treatment-planning pipeline, (2) evaluate GPU-accelerated QIHD on IMRT and IMPT cases with disease-site heterogeneity, and (3) characterize solution times and dose-volume histograms for clinically sized problems.

\section{Methods}
\subsection{Treatment-planning optimization formulation}
Let $V$ be the set of voxels, $V_T\subset V$ the target voxels, $O$ the set of organs at risk (OARs), and $V_o$ the voxels in OAR $o\in O$. Let $B$ be the set of photon control points or proton spots, and $d_{ij}$ denote the dose delivered to voxel $i$ by unit intensity of control point (spot) $j$. The dose-influence coefficients are assembled in $D\in\mathbb{R}^{|V|\times|B|}$. Dose-influence matrices were generated with matRad (Patch 3.2.3) where applicable.\cite{wieser2017development} For nonnegative intensity vector $x$, the dose to voxel $i$ is
$
\mathcal{D}_i=\sum_{j\in B}d_{ij}x_j.
$

For the target, $T_{\min}$ and $T_{\max}$ denote prescribed lower and upper dose bounds and $s_i,s_i'\ge0$ quantify deviations. For OAR $o$, $(UD_o,UV_o)$ denotes a dose--volume requirement in which at most the fraction $UV_o$ of voxels may exceed dose $UD_o$ $(V_{UD_o}\le UV_o)$. With binary indicator $y_i$ and a valid upper dose bound $M$, the planning model is
\begin{align}
\min\quad & \sum_{i\in V_T}(s_i+s_i') \\
\text{s.t.}\quad
&T_{\min}-s_i'\le \sum_{j\in B}d_{ij}x_j\le T_{\max}+s_i, &&i\in V_T,\\
&\sum_{j\in B}d_{ij}x_j\le UD_o+(M-UD_o)y_i, &&i\in V_o,\ o\in O,\\
&\sum_{i\in V_o}y_i\le UV_on_o, &&o\in O,\\
&x,s,s'\ge0,\qquad y_i\in\{0,1\},
\end{align}
where $n_o=|V_o|$. Multiple DVH points can be represented by repeating the OAR constraints. This direct formulation is attractive because the prescription is represented explicitly, but the voxel-level binary variables make realistic instances large MIO problems.\cite{zinchenko2008controlling}

\subsection{End-to-end computational workflow}
\begin{figure}[t]
\centering
\includegraphics[width=0.48\linewidth]{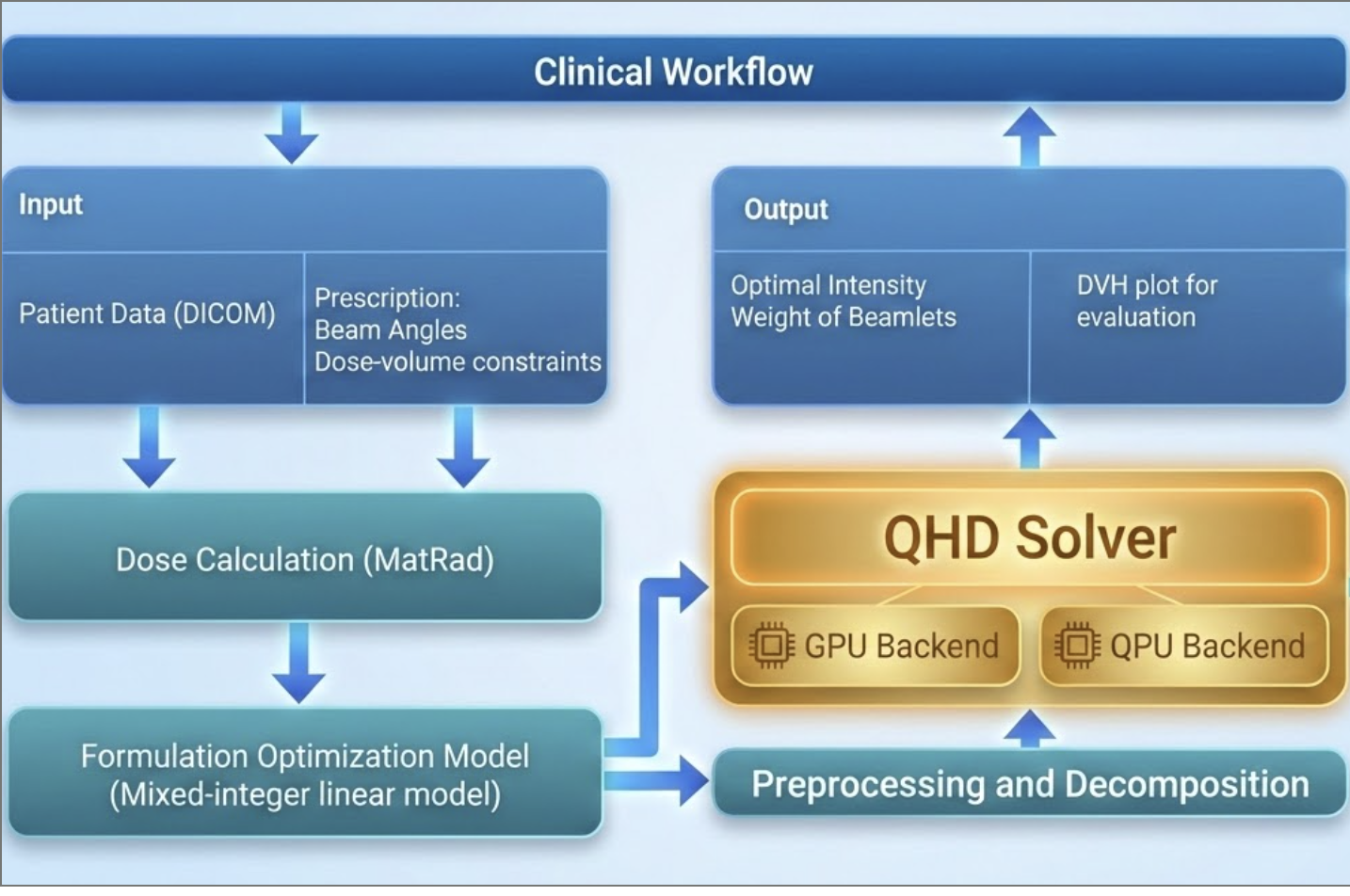}\hfill
\includegraphics[width=0.48\linewidth]{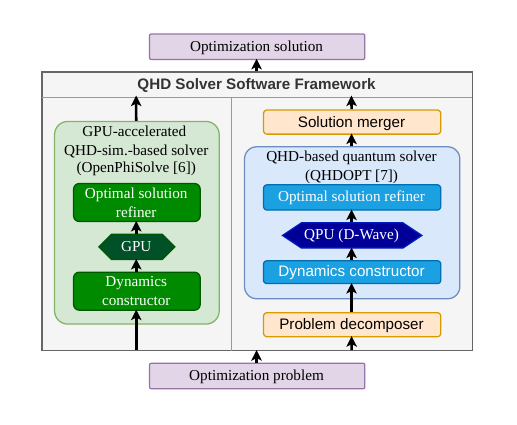}
\caption{End-to-end treatment-planning workflow (left) and modular solver architecture (right). The present clinical-scale experiments use the GPU-accelerated QIHD backend.}
\label{fig:workflow}
\end{figure}

The implemented workflow (Fig.~\ref{fig:workflow}) received imaging-derived dose data, beam configurations, target and organ-at-risk (OAR) structures, and prescription parameters. It then constructed the mathematical formulation, invoked
an optimization backend, and returned optimized control point intensities for photon therapy or spot intensities for proton therapy, along with dose--volume histograms (DVHs) for evaluation. The architecture separated data and model construction from the optimization backend, allowing the same planning model to be evaluated using a CPU-based MIO solver, GPU-accelerated QIHD, or potential future QHD backends. The dominant computational operations involved matrix--vector products and gradient
evaluations using the dose influence matrix. These
operations were implemented in JAX 0.12.0, a Python library for accelerator-oriented array computation, enabling parallel evaluation of trajectories and linear algebra operations on a GPU \cite{chaudhary2025quantum}.

The connection between QIHD and QHD is algorithmic rather than an assertion of equivalent physical computation. QHD evolves a quantum state under a Hamiltonian containing kinetic and objective-dependent potential terms, whereas QIHD evolves classical phase-space trajectories.\cite{leng2023quantum,chaudhary2025quantum} QIHD is used here because it can operate at the scale required for the treatment planning experiments. More details about QIHD can be found in Appendix~1.

\subsection{Datasets and computational evaluation}
The evaluation included four IMRT patient cases (two prostate, one head-and-neck, and one liver) and four IMPT patient cases (two lung, one prostate, and one liver) assembled from the open-source and retrospective de-identified datasets. For each case, the resulting MIO dimensions were recorded. QIHD was run on an NVIDIA GeForce RTX 4060 GPU in a system with an Intel Core Ultra 7 155H CPU. The primary computational endpoints were wall-clock QIHD solution time, number of variables, and the resulting DVHs.

Since the current treatment planning systems are unable to solve these MIO formulations, for a better comparison, the same MIO formulations were also supplied to Gurobi 13.0.3 on the CPU. This benchmark evaluated the difficulty of the explicit formulation of MIO at the voxel-level; it was not intended as a head-to-head benchmark against proprietary optimization algorithms embedded in clinical treatment-planning systems. 

The retrospective datasets used in this study were publicly available and de-identified prior to computational analysis. Therefore, institutional review board (IRB) approval was not required, subject to applicable institutional policies.

\section{Results}
\subsection{IMRT cases}
The four IMRT formulations had 27,621 to 112,565 total variables, with 3,191--100,106 binary variables (Table~\ref{tab:imrt}). QIHD solution times ranged from 5 to 32 min. The largest model by total variable count was the liver case (112,565 variables), whereas the head-and-neck case contained the largest number of the nonzeros elements in dose matrix (70,247,222). These differences illustrated that QIHD was applied across models with substantially different dimensions and sparsity patterns.

\begin{table}[h]
\caption{IMRT model dimensions and QIHD wall-clock solution times.}
\label{tab:imrt}
\scriptsize
\begin{ruledtabular}
\begin{tabular}{lrrrrr}
Case & Continuous & Binary & DVH & Matrix & Time\\
     & Variables & Variables & Constraints & NNZ & (min)\\
\hline
Prostate 1 & 27,621 & 13,360 & 26,902 & 4,570,582 & 8\\
Prostate 2 & 34,421 & 26,488 & 32,426 & 38,163,659 & 5\\
Head-and-neck & 85,651 & 3,191 & 80,547 & 70,247,222 & 32\\
Liver & 112,565 & 100,106 & 111,644 & 20,599,480 & 19\\
\end{tabular}
\end{ruledtabular}
\end{table}

\begin{figure}[h]
\centering
\includegraphics[width=0.44\linewidth]{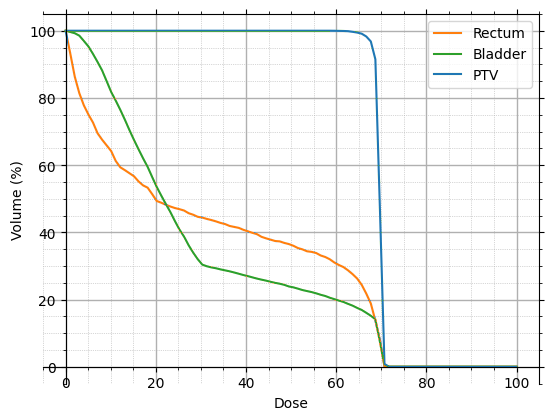}\hfill
\includegraphics[width=0.44\linewidth]{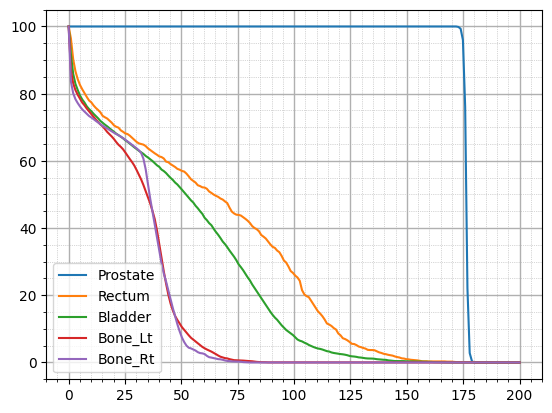}\hfill
\includegraphics[width=0.44\linewidth]{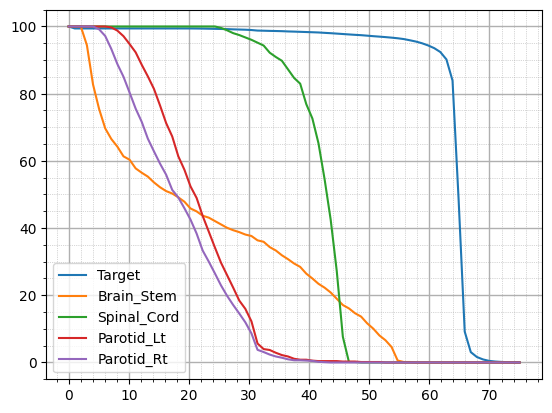}\hfill
\includegraphics[width=0.44\linewidth]{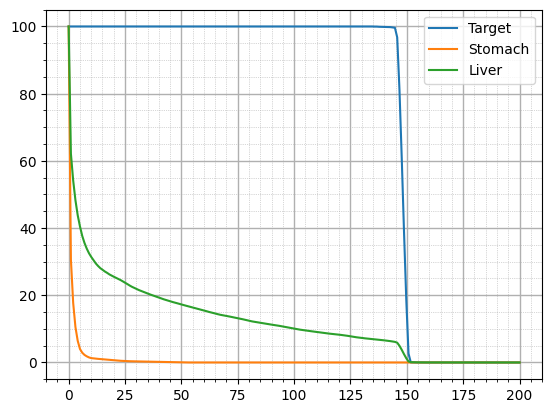}
\caption{Dose--volume histograms generated by QIHD for four IMRT cases: prostate 1, prostate 2, head-and-neck, and liver (left to right).}
\label{fig:imrt}
\end{figure}

Figure~\ref{fig:imrt} shows the DVHs returned by QIHD for the four IMRT cases. Across the cases, the target and OAR curves followed the dose-volume patterns encoded in the corresponding prescriptions, supporting the ability of the pipeline to translate clinical constraints into optimized dose distributions.

\subsection{IMPT cases}
The IMPT formulations contained 69,861--221,342 variables and QIHD solution times ranged from 11 to 35 min (Table~\ref{tab:impt}). The lung 1 case had the largest number of variables and also the largest binary variable count, with 221,342 total variables and 165,715 binary variables. 

\begin{table}[h]
\caption{IMPT model dimensions and QIHD wall-clock solution times.}
\label{tab:impt}
\scriptsize
\setlength{\tabcolsep}{3pt}
\begin{ruledtabular}
\begin{tabular}{lrrrrr}
Case & Continuous & Binary & DVH & Matrix & Time\\
     & Variables & Variables & Constraints & NNZ & (min)\\
\hline
Lung 1   & 221,342 & 165,715 & 200,086 & 4,298,788 & 11\\
Lung 2   & 69,861 & 14,377 & 48,698 & 9,271,187 & 30\\
Prostate & 84,126 & 25,195 & 38,763 & 36,763,332 & 14\\
Liver    & 96,931 & 79,911 & 91,449 & 12,280,133 & 35\\
\end{tabular}
\end{ruledtabular}
\end{table}

\begin{figure}[h]
\centering
\includegraphics[width=0.44\linewidth]{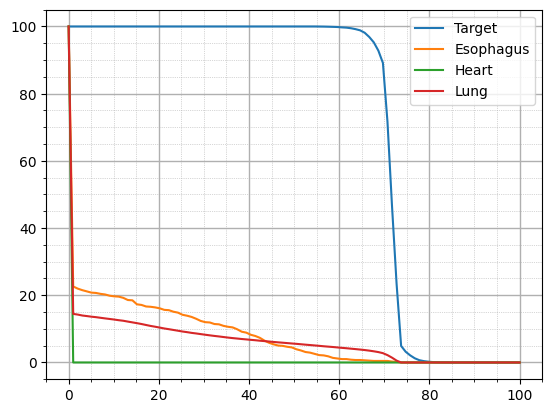}\hfill
\includegraphics[width=0.44\linewidth]{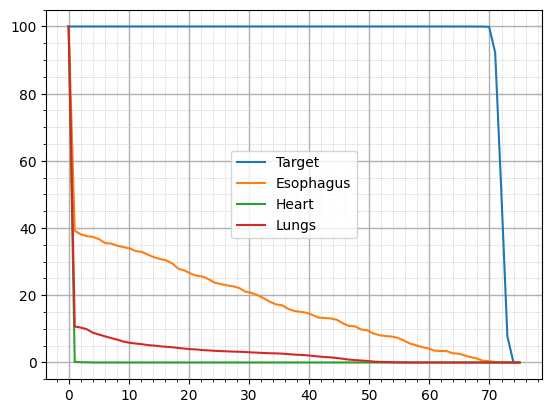}\hfill
\includegraphics[width=0.44\linewidth]{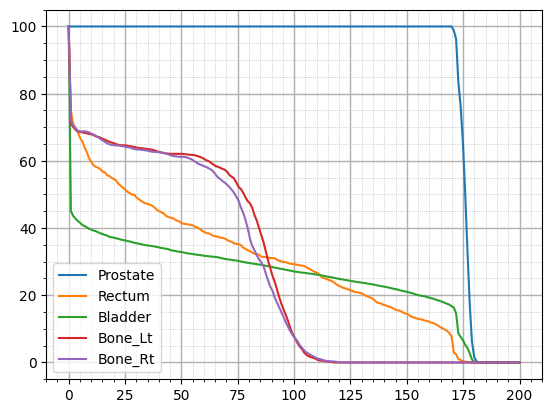}\hfill
\includegraphics[width=0.44\linewidth]{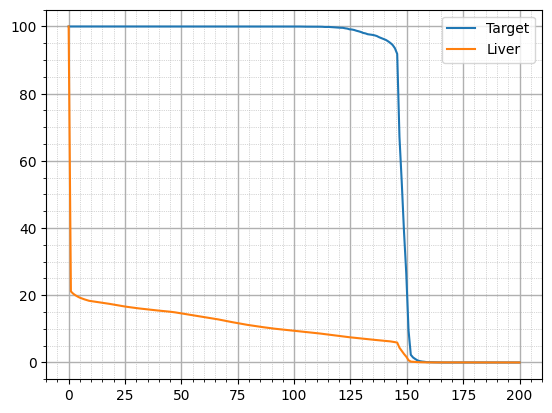}
\caption{Dose--volume histograms generated by QIHD for four IMPT cases: lung 1, lung 2, prostate, and liver (left to right).}
\label{fig:impt}
\end{figure}

\begin{figure}[t]
\centering
\includegraphics[width=0.52\linewidth, clip=true, trim={0 10 0 0}]{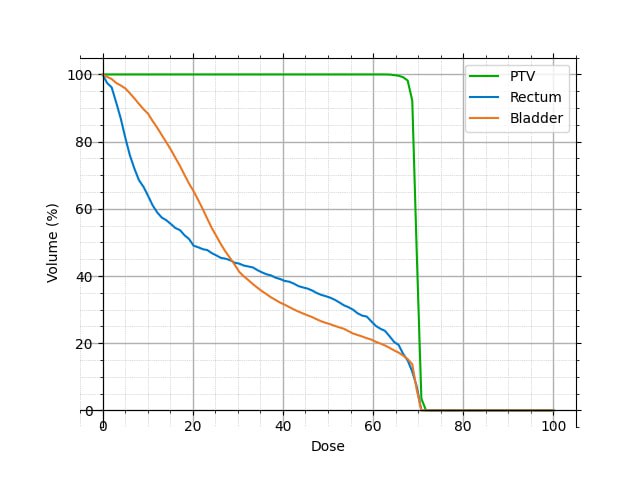}\hfill
\includegraphics[width=0.45\linewidth]{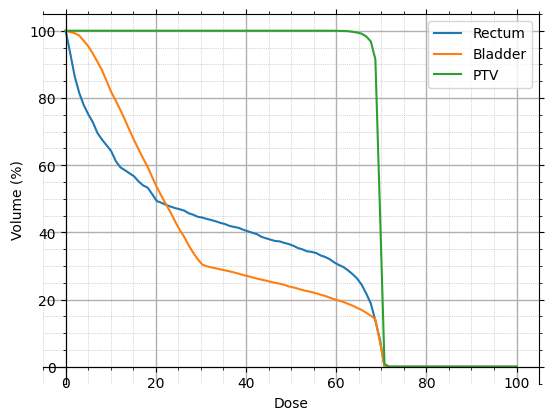}
\caption{Illustrative prostate DVHs from the original NIH Challenge proposal: conventional MIO solution (left) and QIHD solution (right). The reported wall-clock times were approximately 75 and 8 min, respectively.}
\label{fig:legacycomparison}
\end{figure}

\section{Discussion}

The numerical results demonstrate the potential of Quantum-Inspired
Hamiltonian Descent (QIHD) for large-scale radiation therapy treatment
planning. Across eight IMRT and IMPT cases, QIHD generated treatment
plans with dose--volume histograms (DVHs) consistent with the
prescribed clinical requirements. The method solved mixed-integer
optimization (MIO) formulations containing up to 221,342 continuous
variables and 165,715 binary variables within 5--35 minutes on a
single NVIDIA RTX 4060 GPU. These findings demonstrate promising
computational performance using relatively modest hardware resources.

A key advantage of the proposed framework is the direct incorporation
of clinical dose--volume constraints into the optimization model.
Conventional treatment planning approaches often rely on surrogate
objectives and iterative adjustment of optimization parameters,
which may require substantial planner intervention. In contrast,
the proposed formulation explicitly represents dose--volume
requirements using binary variables, potentially reducing manual
tuning and enabling a more automated planning workflow. Despite
the combinatorial complexity of these formulations, QIHD obtained
promising solutions within clinically relevant computational times.
For comparison, the commercial MIO solver Gurobi was unable to
obtain comparable solutions within several hours for most tested
cases. Although this observation suggests potential computational
advantages, systematic benchmarking with solver logs, optimality
gaps, and repeated runs is needed to establish a quantitative
performance comparison.

The successful application of QIHD to both IMRT and IMPT demonstrates
its flexibility across different treatment modalities and anatomical
sites. An important direction for future research is the integration
of robust optimization, particularly for IMPT, where uncertainties
in patient setup and proton range can substantially affect dose
distributions \cite{bortfeld2008robust,unkelbach2018robust}.
Incorporating uncertainty scenarios or ellipsoidal uncertainty
sets into the mixed-integer formulation would allow treatment plans
to account explicitly for these uncertainties. Although this
extension would increase computational complexity, the present
results suggest that QIHD may provide a promising computational
foundation for robust treatment planning.

QIHD is a classical GPU-based algorithm inspired by Quantum
Hamiltonian Descent (QHD), and the present results do not establish
quantum computational advantage. Nevertheless, the connection
to QHD motivates future investigation of quantum and hybrid
quantum-classical implementations. In particular, decomposition
techniques may enable smaller optimization subproblems to be
addressed using quantum hardware while retaining large-scale
linear algebra operations on classical computers.

The present study is limited to eight patient cases and primarily
establishes computational feasibility. Further validation using
standardized dosimetric metrics, comparisons with clinical treatment
planning systems, larger patient cohorts, and treatment deliverability
assessments is needed to establish clinical advantages. Overall,
the results suggest that QIHD offers a promising approach for
efficiently solving large-scale dose--volume-constrained treatment
planning problems, with potential benefits for planning automation
and future extensions to robust optimization.
\section{Conclusions}
A mixed-integer quantum-inspired optimization framework was developed and evaluated on four IMRT and four IMPT cases. GPU-accelerated Quantum-Inspired Hamiltonian Descent (QIHD) was employed to directly incorporate dose--volume constraints into treatment plan optimization. QIHD solved formulations containing up to 221,342 variables within 5--35 minutes on an NVIDIA RTX 4060 GPU, producing dose--volume histograms consistent with the prescribed planning constraints. The resulting treatment plans demonstrated promising dose conformity. Unlike conventional treatment planning approaches that often require iterative adjustment of objective weights and optimization parameters, the proposed framework directly incorporated clinical dose--volume requirements, potentially reducing manual tuning and computational effort. These findings suggest that QIHD offers a promising approach for efficiently solving large-scale mixed-integer treatment planning problems and accelerating the treatment planning workflow. Nevertheless, further evaluation on larger patient cohorts, using standardized dosimetric endpoints, direct comparisons with clinical treatment planning systems, robustness analyses, and quantum hardware experiments, is necessary to establish its clinical advantages and assess the potential for quantum computational speedup.

\section*{Acknowledgments}

This research was conducted in connection with the Quantum Computing Challenge organized by the National Center for Advancing Translational Sciences (NCATS) at the National Institutes of Health (NIH), as part of the Quantum Computing for Biomedical Sciences (Qu-BIT) initiative. The authors' project was awarded as a winner in Phase 1 and Phase 2, Milestone 1, of the challenge. Additionally, Additionally, Tamás Terlaky acknowledges support from Varian Medical Systems, Inc for providing their system.

Generative artificial intelligence tools, including OpenAI's ChatGPT and Anthropic's Claude, were used to assist with language editing and manuscript preparation, as well as code development and debugging. All scientific ideas, mathematical formulations, numerical results, and interpretations were developed or reviewed by the authors. The authors take full responsibility for the accuracy and integrity of the manuscript, including its scientific content, references, computational analyses, and final wording.

\section*{Data Availability Statement}
The treatment planning datasets used in this study were obtained
from two publicly available research resources: PortPy
(Port Planning in Python) and the Common Optimization for
Radiation Therapy (CORT) dataset. 
\cite{craft2014cort,jhanwar2023portpy}

\bibliographystyle{unsrt}
\bibliography{references}

@inproceedings{chaudhary2025quantum,
  title={Quantum-Inspired {Hamiltonian} Descent for Mixed-Integer Quadratic Programming},
  author={Chaudhary, Shreya and Cheng, Jinglei and Kushnir, Samuel and Leng, Jiaqi and Liu, Pengyu and Peng, Yuxiang and Wang, Hanrui and Wu, Xiaodi},
  booktitle={NeurIPS Workshop on GPU-Accelerated and Scalable Optimization},
  year={2025}
}

@article{kushnir2025qhdopt,
  title={{QHDOPT}: A software for nonlinear optimization with quantum {Hamiltonian} descent},
  author={Kushnir, Samuel and Leng, Jiaqi and Peng, Yuxiang and Fan, Lei and Wu, Xiaodi},
  journal={INFORMS Journal on Computing},
  volume={37},
  number={1},
  pages={107--124},
  year={2025},
  publisher={INFORMS}
}

@article{craft2012multicriteria,
  title={Multicriteria {VMAT} optimization},
  author={Craft, David and McQuaid, Dualta and Wala, Jeremiah and Chen, Wei and Salari, Ehsan and Bortfeld, Thomas},
  journal={Medical Physics},
  volume={39},
  number={2},
  pages={686--696},
  year={2012},
  publisher={Wiley Online Library}
}

@article{wieser2017development,
  title={Development of the open-source dose calculation and optimization toolkit MatRad},
  author={Wieser, Hans-Peter and Cisternas, Eduardo and Wahl, Niklas and Ulrich, Silke and Stadler, Alexander and Mescher, Henning and M{\"u}ller, Lucas-Raphael and Klinge, Thomas and Gabrys, Hubert and Burigo, Lucas},
  journal={Medical Physics},
  volume={44},
  number={6},
  pages={2556--2568},
  year={2017},
  publisher={Wiley Online Library}
}

@article{zinchenko2008controlling,
  title={Controlling the dose distribution with {gEUD}-type constraints within the convex radiotherapy optimization framework},
  author={Zinchenko, Yuriy and Craig, Timothy and Keller, Harald and Terlaky, Tam{\'a}s and Sharpe, Michael B},
  journal={Physics in Medicine \& Biology},
  volume={53},
  number={12},
  pages={3231},
  year={2008},
  publisher={IOP Publishing}
}

@article{leng2023quantum,
  title={Quantum {Hamiltonian} descent},
  author={Leng, Jiaqi and Hickman, Ethan and Li, Joseph and Wu, Xiaodi},
  journal={arXiv preprint arXiv:2303.01471},
  year={2023}
}

@article{heijmen2018fully,
  title={Fully automated, multi-criterial planning for Volumetric Modulated Arc Therapy--An international multi-center validation for prostate cancer},
  author={Heijmen, Ben and Voet, Peter and Fransen, Dennie and Penninkhof, Joan and Milder, Maaike and Akhiat, Hafid and Bonomo, Pierluigi and Casati, Marta and Georg, Dietmar and Goldner, Gregor and others},
  journal={Radiotherapy and Oncology},
  volume={128},
  number={2},
  pages={343--348},
  year={2018},
  publisher={Elsevier}
}

@article{hussein2018automation,
  title={Automation in intensity modulated radiotherapy treatment planning—a review of recent innovations},
  author={Hussein, Mohammad and Heijmen, Ben JM and Verellen, Dirk and Nisbet, Andrew},
  journal={The British Journal of Radiology},
  volume={91},
  number={1092},
  pages={20180270},
  year={2018},
  publisher={The British Institute of Radiology.}
}

@article{mcintosh2017fully,
  title={Fully automated treatment planning for head and neck radiotherapy using a voxel-based dose prediction and dose mimicking method},
  author={McIntosh, Chris and Welch, Mattea and McNiven, Andrea and Jaffray, David A and Purdie, Thomas G},
  journal={Physics in Medicine \& Biology},
  volume={62},
  number={15},
  pages={5926--3944},
  year={2017},
  publisher={IOP Publishing}
}

@article{breedveld2012icycle,
  title={{iCycle:} Integrated, multicriterial beam angle, and profile optimization for generation of coplanar and noncoplanar {IMRT} plans},
  author={Breedveld, Sebastiaan and Storchi, Pascal RM and Voet, Peter WJ and Heijmen, Ben JM},
  journal={Medical Physics},
  volume={39},
  number={2},
  pages={951--963},
  year={2012},
  publisher={Wiley Online Library}
}

@article{unkelbach2018robust,
  title={Robust radiotherapy planning},
  author={Unkelbach, Jan and Alber, Markus and Bangert, Mark and Bokrantz, Rasmus and Chan, Timothy CY and Deasy, Joseph O and Fredriksson, Albin and Gorissen, Bram L and Van Herk, Marcel and Liu, Wei and others},
  journal={Physics in Medicine \& Biology},
  volume={63},
  number={22},
  pages={22TR02},
  year={2018},
  publisher={IOP Publishing}
}

@article{bortfeld2008robust,
  title={Robust management of motion uncertainty in intensity-modulated radiation therapy},
  author={Bortfeld, Thomas and Chan, Timothy CY and Trofimov, Alexei and Tsitsiklis, John N},
  journal={Operations Research},
  volume={56},
  number={6},
  pages={1461--1473},
  year={2008},
  publisher={INFORMS}
}

@article{shepard1999optimizing,
  title={Optimizing the delivery of radiation therapy to cancer patients},
  author={Shepard, David M and Ferris, Michael C and Olivera, Gustavo H and Mackie, T Rockwell},
  journal={{SIAM} Review},
  volume={41},
  number={4},
  pages={721--744},
  year={1999},
  publisher={{SIAM}}
}

@inproceedings{mohammadisiahroudi2024optimization,
  title={On Optimization Challenges in Proton Therapy Treatment Planning},
  author={Mohammadisiahroudi, Mohammadhossein and Zou, Jennifer Wei and Dong, Lei and Zinchenko, Yuriy and Terlaky, Tam{\'a}s},
  booktitle={AAPM 66th Annual Meeting \& Exhibition},
  year={2024},
  organization={AAPM}
}

@article{rahimi2026quantum,
  title={Quantum computing for radiation therapy optimization},
  author={Rahimi, Robabeh and SaiToh, Akira and Modiri, Arezoo and Nakano, Yuichiro and Okada, Ken N and Tsukano, Satoyuki and Zhang, Baoshe and Fujii, Keisuke and Kitagawa, Masahiro and Sawant, Amit},
  journal={Medical Physics},
  volume={53},
  number={1},
  pages={e70269},
  year={2026},
  publisher={Wiley Online Library}
}

@article{pakela2020quantum,
  title={Quantum-inspired algorithm for radiotherapy planning optimization},
  author={Pakela, Julia M and Tseng, Huan-Hsin and Matuszak, Martha M and Ten Haken, Randall K and McShan, Daniel L and El Naqa, Issam},
  journal={Medical Physics},
  volume={47},
  number={1},
  pages={5--18},
  year={2020},
  publisher={Wiley Online Library}
}

@article{jhanwar2023portpy,
  title={PortPy: An Open-Source Python Package for Planning
         and Optimization in Radiation Therapy Including
         Benchmark Data and Algorithms},
  author={Jhanwar, Gourav and Tefagh, Mojtaba and
          Taasti, Vicki T. and Alam, Sadegh R. and
          Tuomaala, Seppo and Nadeem, Saad and
          Zarepisheh, Masoud},
  journal={AAPM 65th Annual Meeting \& Exhibition},
  year={2023}
}

@article{craft2014cort,
  title={Shared data for intensity modulated radiation therapy
         ({IMRT}) optimization research: the {CORT} dataset},
  author={Craft, David and Bangert, Mark and Long, Troy
          and Papp, D{\'a}vid and Unkelbach, Jan},
  journal={GigaScience},
  volume={3},
  number={1},
  pages={37},
  year={2014},
  doi={10.1186/2047-217X-3-37}
}

@techreport{mohammadisiahroudiadaptive,
  author      = {Mohammadisiahroudi, Mohammadhossein
                 and Sampourmahani, Pouya
                 and Zou, Wei
                 and Dong, Lei
                 and Zinchenko, Yuriy
                 and Terlaky, Tam{\'a}s},
  title       = {An Adaptive Relaxation Technique for Mixed-integer
                 Formulation of Proton Therapy Treatment Optimization},
  institution = {Department of Industrial and Systems Engineering,
                 Lehigh University},
  type        = {ISE Technical Report},
  number      = {25T-008},
  year        = {2025},

  url         = {https://engineering.lehigh.edu/sites/engineering.lehigh.edu/files/_DEPARTMENTS/ise/pdf/tech-papers/25/25T_008.pdf}
}
\newpage
\section{Quantum-Inspired Hamiltonian Descent}
\label{app:QIHD}
QIHD simulates a classical dynamics governed by time-dependent Hamiltonian
\begin{equation}
H(x,y,t)=e^{\varphi(t)}\left(\frac{1}{2}\|y\|_2^2\right)+e^{\chi(t)}f(x),
\label{eq:qihd}
\end{equation}
where $x$ is the position (optimization variables), $y$ is the conjugate momentum, and $f(x)$ is the penalized treatment-planning objective. The dynamics obey Hamilton's equations,
\begin{equation}
\dot{x}=\nabla_yH,\qquad \dot{y}=-\nabla_xH.
\end{equation}
Randomized initial positions and momenta provide multiple candidate trajectories, analogous to sampling multiple candidate solutions. A first-order discretization can be written as
\begin{align}
y_k &= y_{k-1}-\nabla_xH(x_{k-1},y_{k-1},t_{k-1})\Delta t,\\
x_k &= x_{k-1}+\nabla_yH(x_{k-1},y_{k-1},t_{k-1})\Delta t.
\end{align}
\section{Treatment Planning Optimization Parameters}
\label{app:clinical_parameters}

Table~\ref{tab:clinical_parameters} summarizes the target
dose bounds and organ-at-risk (OAR) dose--volume constraints
used in the optimization formulations. The parameter values
were extracted from the computational implementations for
the prostate, head-and-neck, lung, and liver cases.

For each OAR, $V_d$ denotes the percentage of voxels
receiving a dose of at least $d$. The upper dose--volume
constraint $V_d \leq v$ requires that no more than $v\%$
of the corresponding organ volume receive dose $d$
or greater.

\begin{table*}[t]
\caption{Target dose bounds and OAR dose--volume constraints
specified in the optimization models. Dose values are
reported in the numerical units used in the source code;
conversion to physical dose units requires verification.}
\label{tab:clinical_parameters}
\footnotesize
\setlength{\tabcolsep}{5pt}
\begin{ruledtabular}
\begin{tabular}{lllll}
Disease site & Structure & Lower bound &
Upper bound & Dose--volume constraint \\
\hline
Liver & Target & 145 & 150 & -- \\
      & Stomach & -- & -- & $V_{20}\leq 10\%$ \\
      & Liver & -- & -- & $V_{80}\leq 20\%$ \\
\hline
Lung 1 and 2 & Target & 70 & 72 & -- \\
      & Esophagus & -- & -- & $V_{60}\leq 17\%$ \\
      & Heart & -- & -- & $V_{50}\leq 25\%$ \\
      & Lungs & -- & -- & $V_{50}\leq 10\%$ \\
\hline
Head and neck & Target & 63 & 65 & -- \\
      & Brainstem & -- & -- & $V_{54}\leq 5\%$ \\
      & Spinal cord & -- & -- & $V_{45}\leq 5\%$ \\
      & Left parotid & -- & -- & $V_{30}\leq 50\%$ \\
      & Right parotid & -- & -- & $V_{30}\leq 50\%$ \\
\hline
Prostate & Target & 171 & 180 & -- \\
      & Rectum & -- & -- & $V_{170}\leq 20\%$ \\
      & Bladder & -- & -- & $V_{171}\leq 35\%$ \\
      & Left femoral bone & -- & -- & $V_{122}\leq 3\%$ \\
      & Right femoral bone & -- & -- & $V_{122}\leq 3\%$ \\
\end{tabular}
\end{ruledtabular}
\end{table*}

The target dose bounds were incorporated into the
optimization objective through nonnegative deviation
variables, while the OAR dose--volume requirements
were represented explicitly using binary variables.
The Big-M parameters were set to 200 for liver and
prostate, 72 for lung, and 70 for head-and-neck cases.

\end{document}